\documentclass[microtype]{jloganal}

\usepackage{amsmath,amssymb,enumerate,bbm}

\newcommand{\assign}{:=}
\newcommand{\mathd}{\mathrm{d}}
\newcommand{\tmem}[1]{{\em #1\/}}
\newcommand{\tmop}[1]{\ensuremath{\operatorname{#1}}}
\newcommand{\tmstrong}[1]{\textbf{#1}}
\newcommand{\tmtextit}[1]{{\itshape{#1}}}
\newcommand{\um}{-}
\newcommand{\upl}{+}
\newenvironment{descriptioncompact}{\begin{description} }{\end{description}}
\newenvironment{enumeratenumeric}{\begin{enumerate}[1.] }{\end{enumerate}}
\newtheorem{definition}{Definition}
\newtheorem{lemma}{Lemma}
\newtheorem{proposition}{Proposition}
\newtheorem{theorem}{Theorem}

\newcommand{\compr}[2]{\ensuremath{\left\{ #1  \mid  #2 \right\}}}
\newcommand{\Max}[1]{\ensuremath{\tmop{Max} ( #1)}}
\newcommand{\N}{\ensuremath{\mathbbm{N}}}

\newcommand{\integral}[1]{\ensuremath{P ( #1)}}
\newcommand{\intl}[2]{\ensuremath{[ #1 < I ( #2)]}}

\newcommand{\Top}{1}
\newcommand{\Bot}{0}

\begin{document}

\title{Integrals and valuations}
\arxivreference{math.LO/0808.1522}

\author{Thierry Coquand}
\givenname{Thierry }
\surname{Coquand}
\address{Computing Science Department, G\"oteborg University}
\email{coquand@chalmers.se}
\urladdr{}

\author{Bas Spitters}
\givenname{Bas}
\surname{Spitters}
\address{Mathematics Department, Eindhoven University of Technology}
\email{spitters@cs.ru.nl}
\urladdr{}

\begin{abstract}
  We construct a homeomorphism between the compact regular locale of integrals
  on a Riesz space and the locale of measures(valuations) on its spectrum. In
  fact, we construct two geometric theories and show that they are
  biinterpretable. The constructions are elementary and tightly connected to
  the Riesz space structure.
\end{abstract}
\maketitle

\section{Introduction}

The goal of this paper is to give a constructive formulation of the Riesz
representation theorem. The Riesz representation theorem states that there is
an isomorphism between {\tmem{integrals}} and {\tmem{regular measures}} on
compact spaces. An integral on a compact Hausdorff space $X$ is a positive
linear functional $I : C (X) \rightarrow \R$ (and we shall consider only maps
such that $I (1) = 1$). A regular measure can be identified with a
{\tmem{continuous valuation}}, where a {\tmem{valuation}} on $X$ is a map $\mu
: O (X) \rightarrow [0, 1]$ which is monotone, if $U \subseteq V$ then $\mu
(U) \leqslant \mu (V)$, and such that $\mu (\emptyset) = 0$ and $\mu (U \cap
V) + \mu (U \cup V) = \mu (U) + \mu (V)$ and $\mu (X) = 1$. The continuity
condition demands that $\mu (U)$ is the sup of $\mu (V)$ for $V$ well-inside
$U$ (i.e. such that $U$ contains the closure of $V$). An equivalent way to
express this condition is to state the \tmtextit{continuity} property: if
$V_i$ is a directed family then $\mu ( \bigcup V_i) = \sup \mu (V_i)$. A
subset is {\tmem{directed}} if it is inhabited and every two elements have a
common upper bound. Such continuous valuations extend uniquely to Borel
measures; see~{\cite{Jung}} for an overview.

From a constructive point of view there is a crucial difference between the
two notions. We will now outline these differences; precise definitions can be
found below. The integral $I (f)$ of a function $f \in C (X)$ is a
\tmtextit{Dedekind real}. Intuitively, this means that one can compute
arbitrary rational approximations. This may not be the case for the valuation
$\mu (U)$ of an open $U$: in general we do not have the property that for $r <
s$,
\[ \mu (U) < s \vee r < \mu (U) . \]
Constructively the valuation $\mu (U)$ is only a \tmtextit{lower real}, and
can be thought of as a predicate $r < \mu (U)$ on the rationals. This
predicate is downward closed: if $r < \mu (U)$ and $s \leq r$ then we have $s
< \mu (U)$, but in general, given $\epsilon > 0$ we are not given a way to
compute a rational $\epsilon$ approximation of $\mu (U)$. Given an integral
$I$ we define a corresponding valuation $\mu_I (U)$ by taking the sup of $I
(f)$ over all $0 \leqslant f \leqslant 1$ the support of which is included in
$U$. It is remarkable that for \tmtextit{any} valuation $\mu$ one can
conversely find a (unique) integral $I$ such that $\mu = \mu_I$. So despite
the fact that one may not be able to compute $\mu (U)$, it is still possible
to compute $\int fd \mu$ as a {\tmem{Dedekind real}} as the supremum of
\[ \sum s_i \mu (s_i < f < s_{i + 1}) \]
over all partitions $s_0 < \cdots < s_n$ of the range $f ([a, b])$. A priori
this supremum will only be a {\tmem{lower real}}.

As usual in constructive mathematics all structures carry a natural, but
implicit, topology and all constructions are continuous. To make this
structure explicit we start from a Riesz space $R$ and associate three formal
spaces to it that are all compact completely regular: the maximal spectrum
$\Max{R} = X$ (intuitively, $R$ is then a dense subset of $C (X)$), the space
of integrals $\tmop{INT} (R)$ and the space of valuations $\tmop{VAL} (X)$.
All three spaces are defined as \tmtextit{propositional geometrical theories}.
A geometric formula is one of the form $\psi \Rightarrow \varphi$, where the
formulas $\psi$ and $\varphi$ are positive, {\tmem{i.e.}} they are built up
from atomic formulas using only (finite) conjunction, (infinite) disjunction.
A {\tmem{geometric theory}} is a theory all of whose axioms are geometric. The
main point of this paper is to define two interpretability maps, showing how
to interpret the theory $\tmop{VAL} (X)$ in the theory $\tmop{INT} (R)$
(intuitively how to define the measure from an integral) and how to interpret
the theory $\tmop{INT} (R)$ in the theory $\tmop{VAL} (X)$ (intuitively how to
define the measure from the integral). The Riesz representation theorem can
then be stated as the fact that these two maps define an homeomorphism between
the corresponding formal spaces $\tmop{VAL} (X)$ and $\tmop{INT} (R)$ where
the topology on $\tmop{INT} (R)$ is the weak topology. Hence we arrive at a
concrete constructive statement of the Riesz representation theorem which is
valid in any topos.

The present article is part of our program to apply the logical approach to
abstract algebra~{\cite{Coquand/Lombardi}} to (functional)
analysis~{\cite{Coquand:Stone,Coquand/Spitters:StoneY,Spitters:caitnoac,Coquand:HahnBanach}}.
It may be seen as a contribution to Hilbert's program of logically translating
the use of infinitary methods to finitary, or constructive, ones. It is also
continuation of a tradition in topos theory,
e.g.~{\cite{banaschewskimulvey06}}, but in a more explicit manner.{\footnote{
We avoid the axiom of (countable) choice, and, moreover, we refrain from using
the power set axiom. One may wonder how we treat the set of all real numbers
in such a framework. In fact, we do not use this set at all. We only consider
the {\tmem{formal space}} of real numbers.}} It turns out that our program
sometimes gives shorter proofs of more general results than a direct
constructive treatment in the sense of Bishop. Moreover, the space of
valuations does not naturally carry a metric structure and hence the
topological structure, explicit in our presentation, is hidden in Bishop's
treatment of the Riesz representation theorem. We emphasize, however, that all
our results \tmtextit{are} acceptable by Bishop's standard. Since we do not
construct points, we have no use for the axiom of choice, not even the
countable version which is available to Bishop.

\subsection{Formal measure and integration theory}As outlined
in~{\cite{Coquand:measures-povec,Spitters:caitnoac,Coquand/Palmgren}} a formal
theory of measure and integration may be developed along the following lines.

In a usual set-theoretic foundation of measure theory one considers certain
functions which are defined to be `measurable'. Then relative to a measure one
identifies all the functions which are equal almost everywhere and obtains a
vector lattice $L_0$ of measurable `functions'. Instead, one may consider such
a vector lattice from the beginning, abstracting from the set-theoretic
foundations. The benefits of this approach have been emphasized by Kolmogorov,
Caratheodory and von Neumann~{\cite{Rota}}. In the present article we focus on
the theory of integrals defined on formal functions and valuations defined on
formal {\tmem{opens}}. For a formal treatment of Borel sets we refer
to~{\cite{Coquand/Palmgren,Coquand:measures-povec,Simpson}}.

The abstract space of functions is captured by a Riesz space (a vector
lattice) which we require to have a strong unit{\footnote{An even weaker
requirement would have been to demand that we are given a lattice ordered
Abelian group. Such a group can be extended to a Riesz space over the
rationals; e.g.~{\cite{Coquand:Stone}}.}}. An integral is a
{\tmem{continuous}} linear functional on the Riesz space. On the other hand, a
measure is typically only{\tmem{ lower semi-continuous}}. This suggests that
an integral will be a map to the {\tmem{Dedekind}} reals, but that a valuation
will map to the {\tmem{lower}} reals. The Riesz representation theorem will be
presented in the form of a homeomorphism between the formal space of integrals
on a Riesz space and valuations on the opens of its spectrum. By the
Stone-Yosida theorem any Riesz space $R$ with strong unit can be densely
embedded in the space of continuous functions over its spectrum $\Max{R}$.
This can be proved constructively~{\cite{Coquand/Spitters:StoneY}}. The
integral extends to this space of continuous functions. In this sense our
approach is close to the Daniell integral.

Alternatively, we could have started from a compact completely regular locale
$X$ and construct the Riesz space $C (X)$ of continuous functions on this
locale. Then $\Max{C (X)} \cong X$. However, we also want to include
`syntactic' Riesz spaces such as the Riesz space of rational piecewise linear
functions on $[0, 1]$.

\subsection{Overview}Section~\ref{sec:Riesz} contains the statement of the
Riesz representation theorem, the main result of the article. The statement is
geometric with joins restricted to countable sets. This allows us to use
logical methods to conclude classically that there has to be a constructive
proof. We construct such a proof in Section~\ref{Sec:Proof}. The proof uses a
concrete theory of non-increasing functions, which we call $\Delta$-functions
in Section~\ref{sec:Delta}.

\subsection{Notation}We use the letters $a, b, f, g, ...$ for elements of $R$,
and the letters $x, y, z, ...$ for elements of the lattice $\tmop{Spec} (R)$.
We write {\Top} for the top element of a lattice and {\Bot} for the bottom
element.

\section{Preliminaries}

\subsection{Various kinds of real numbers\label{reals}}

We recall some facts about the real numbers; see
e.g.~{\cite[D4.7]{johnstone02b}}.

\begin{definition}
  \label{def:reals}A {\tmem{lower real}} is an inhabited, down-closed, open
  subset of the rationals. The collection of lower reals is denoted by
  $\R_{\tmop{low}}$. Upper reals are defined similarly and denoted
  $\R_{\tmop{up}}$. An {\tmem{interval}} consists of a pair $(L, U)$ \ of an
  upper real and a lower real such that $L \leqslant U$: if $s$ in $L$ and $t$
  in U, then $s < t$. A {\tmem{Dedekind}} real is an interval $(L, U)$ which
  is arbitrarily small: for every $s < t$, either $s \in L$ or $t \in U$. The
  Dedekind reals will be denoted by $\R$.
\end{definition}

The lower reals include $\upl \infty$. In order to exclude it, we would need
to pose a non-geometric restriction. This issue will not be important in the
rest of the paper.

Lower (likewise upper) reals are closed under addition and closed under
multiplication by a positive rational. The lower and upper reals are
{\tmem{not}} closed under subtraction, but one {\tmem{can}} subtract a lower
real from an upper real and obtain an upper real. The non-strict inequality
$\leqslant$ is given by inclusion of subsets. The supremum of an inhabited set
of lower reals is a lower real. The infimum of an inhabited set of upper reals
is an upper real. This is an important motivation for the use of lower (upper)
reals: the supremum of a sequence of rationals need not be a Dedekind real,
but it {\tmem{is}} a lower real.

In the absence of the powerset operator, the lower reals are better considered
as a formal space rather than a set, but we will not emphasize this point.

In the presence of dependent choice the Dedekind reals coincide with the
Cauchy reals.

\begin{lemma}
  \label{low-up}Let $L$ be a lower real and $U$ be an upper real.

  If $L \leqslant U$, then for all rational $p$, $L + p - U \leqslant p$.

  Conversely, if $L + p - U \leqslant p$ for some rational $p$, then $L
  \leqslant U$.
\end{lemma}

The following lemma will be used a number of times below:

\begin{lemma}
  \label{less-ineq}The relation $a \leqslant b + c$ for lower reals can be
  stated geometrically in two equivalent ways:
  \begin{enumerate}
    \item $p < a \rightarrow \bigvee_{r + s = p} (r < b \wedge s < c)$;

    \item $r + s < a \rightarrow r < b \vee s < c$.
  \end{enumerate}
  A similar statement holds for upper reals.
\end{lemma}

\begin{proof}
  The implication from 1 to 2 is direct. For the implication from 2 to 1 we
  observe that if $p < a$, then there exists $\varepsilon > 0$ such that $p +
  \varepsilon < a$. Choose a rational $q < b$ and a natural number $N$ such
  that $p - q - N \varepsilon < c$. By hypothesis 2, with premiss $p +
  \varepsilon < a$, we have for every $n$, $q + n \varepsilon < b \vee p +
  \varepsilon - (q + n \varepsilon) < c$. Since the first disjunct holds for
  $n = 0$ and the second for $n \geqslant N + 1$, there exists $n$ such that
  $q + n \varepsilon < b$ and $p + \varepsilon - (q + (n + 1) \varepsilon) <
  c$. We can now take $r = q + n \varepsilon$ and $s = p - q + n \varepsilon$.
\end{proof}

The inequality between a lower real and a Dedekind real can also be stated
geometrically.

\subsection{Logic and topology}

In set theory, {\tmem{i.e.}} in the topos Set, one uses topological spaces to
deal with continuity. However, statements including points of topological
spaces are often difficult to generalize to arbitrary toposes. Fortunately, it
is often possible to resort to the lattice structure of the open sets of a
topological space. These complete distributive lattices are thus called
`pointfree' spaces, or locales~(see {\cite{johnstone82}}). In the topos Set
one can often reconstruct the points from this lattice; to be precise, there
is an adjunction between the category of topological spaces and the category
of locales, which restricts to an equivalence of categories between compact
Hausdorff spaces and compact completely regular locales. In general, this
equivalence is not present in a topos. When generalizing theorems from the
topos Set to an arbitrary topos focusing on locales is often the better
choice. One reason for this is that a locale may be defined by geometric
theory. In logical terms the locale is its syntactic category, often called
the Tarski-Lindenbaum algebra {\emdash} that is, the poset of provable
equivalence classes, ordered by provable entailment. The correspondence
between the locale and the theory is the usual completeness and consistency
link between theories and models. The models of the theory correspond to
completely prime filters, {\tmem{i.e.}} points of the locale presented by the
lattice. In this way, a point $x$ in a topological space defines a model of
the corresponding theory: a basic proposition $I$ is true in the model $x$ iff
$x \in I$. This view leads us to consider theories as primary objects of
study; their models, the points, will be derived concepts. Hence topology is
propositional geometric logic; see e.g.~{\cite{johnstone02b,Vic:LocTopSp}}.

\subsection{Spectrum of a Riesz space\label{spectrum}}

\begin{definition}
  An {\tmem{ordered vector space}} is a vector space with a partial order
  $\leqslant$ such that
  \begin{enumerate}
    \item If $x \leqslant y$, then for all $z$, $x + z \leqslant y + z$;

    \item if $0 \leqslant x$, then for all $a \geqslant 0$, $0 \leqslant a x$.
  \end{enumerate}
  A {\tmem{Riesz space{\tmstrong{}}}} {\tmem{(or vector lattice)}} is an
  ordered vector space where the order structure is a lattice. An element $1$
  is a {\tmem{strong unit}} if for all $x$ there exists $n$ such that $- n 1
  \leqslant x \leqslant n 1$.
\end{definition}

As noted in the introduction Riesz spaces provide an algebraic way to talk
continuous functions on a compact completely regular locale.

We will consider Riesz spaces as $\Q$-vector spaces.

In a Riesz space one defines $f^+ \assign f \vee 0$, $f^- \assign 0 \vee (\um
f)$ and $|f| \assign f^+ + f^-$ and derives that $f = f^+ - f^-$.

The spectrum of a Riesz space $R$ is the space of all its representations
{\emdash} Riesz morphisms from $R$ to $\R$. It may be
presented~{\cite{Coquand/Spitters:StoneY}} as the locale which is freely
generated by the collection of tokens $D (a)$, one for each $a$ in $R$,
subject to the following relations:
\begin{enumeratenumeric}
  \item $D (1) = \Top$;

  \item $D (a) \wedge D (- a) = \Bot$;

  \item $D (a + b) \leqslant D (a) \vee D (b)$;

  \item $D (a) = \Bot$, if $a \leqslant 0$;

  \item $D (a \vee b) = D (a) \vee D (b)$;

  \item $D (a) = \bigvee_{s > 0} D (a - s)$.
\end{enumeratenumeric}
As proved in Proposition 2.6 in~{\cite{Coquand:Stone}} one can derive the
relations $D (a) = D (a^+)$ and $D (a \wedge b) = D (a) \wedge D (b)$ from the
ones above. In fact, either of them is equivalent to the relation 5 given
(1-4).

One proves that this locale, $\Max{R}$, is compact and completely regular by
interpreting the geometric theory (1-6) in the coherent theory (1-5) by
interpreting $D (a)$ in (1-6) as $\bigvee_{s > 0} D (a - s)$ in (1-5);
see~{\cite{Coquand:Stone}}. In terms of locales this means that the locale is
a retraction of the coherent locale generated by (1-5). The relations for
$\vee$ and $\wedge$ allow us to reduce $\vee \wedge D (a_{i j})$ to $D (\vee
\wedge a_{i j})$. So the collection of $D (a)$ actually forms a basis, rather
than only a subbasis, for the locale. We write $\tmop{Spec} (R)$ for the
distributive lattice generated by (1-5).

\begin{theorem}
  {\tmem{{\cite{Coquand:Stone}}}} The order in $\tmop{Spec} (R)$ is $D (a)
  \leqslant D (b)$ iff there exists $n$ such that $a^+ \leqslant n b^+$. The
  order in the locale $\Max{R}$ is $D (a) \leqslant D (b)$ iff for all
  $\varepsilon > 0$ there exists $n$ such that $(a - \varepsilon)^+ \leqslant
  n b^+$.
\end{theorem}

The order on $\tmop{Spec} (R)$, as opposed to the order on $\Max{R}$, is
defined geometrically from the order on $R$.

Intuitively, the open $D (a)$ in the locale corresponds to the set
$\compr{\sigma}{\hat{a} (\sigma) > 0}$, where $\hat{a} : \Max{R} \rightarrow
\R$ is the function defined by $\hat{a} (\sigma) \assign \sigma (a)$ for
$\sigma$ in the locale. In the presence of the full axiom of choice this can
be made precise as it allows us to prove that the spectrum has {\tmem{enough}}
points.

\begin{proposition}
  A model $m$ of the geometric theory above, i.e.~a point of the spectrum as a
  locale, defines a representation
  \[ \sigma_m (a) \assign ( \compr{r}{m \models D (a - r)}, \compr{s}{m
     \models D (s - a)}) . \]
\end{proposition}

\begin{proof}
  Lemma 1 in~{\cite{Coquand/Spitters:StoneY}} proves that this defines a
  Dedekind cut. By axiom 5 $\sigma (a \vee b) = \sigma (a) \vee \sigma (b)$.
  By axioms 2,4 $D (1 - (1 - \frac{1}{n})) = \Top$. It follows that $\sigma
  (1) = 1$. As observed in~{\cite{Coquand/Spitters:StoneY}} a map satisfying
  these properties is a representation.
\end{proof}

The Stone-Yosida representation theorem states that there is a embedding of
$R$ into the locale of (Dedekind) real valued continuous functions on its
spectrum which is dense with respect to the sup-norm. The sup-norm is the
upper real $\|a\|$ defined by $\|a\|< \lambda$ iff there exists $\lambda' <
\lambda$ such that $|a| \leqslant \lambda' 1$. A constructive proof of this
theorem can be found in~{\cite{Coquand/Spitters:StoneY}}.{\hspace*{\fill}}

\section{\label{sec:Riesz}Statement of the Riesz-representation theorem}

The goal of this section is to state, in Subsection~\ref{statement}, the Riesz
representation theorem as the existence of a homeomorphism between the formal
compact completely regular spaces of integrals and valuations.
Theorem~\ref{thm:Riesz} contains the proof of the representation theorem.

\subsection{The space of integrals}

Let $R$ be a Riesz space with strong unit 1.

\begin{definition}
  A {\tmem{(probability) integral}} $I$ on a Riesz space $R$ is a positive
  linear functional {\emdash} that is, it is a linear map to the Dedekind
  reals and if $x \geqslant 0$, then $I (x) \geqslant 0$ {\emdash} and such
  that $I (1) = 1$.
\end{definition}

An integral is continuous with respect to the sup-norm: if $|f| \leqslant r$,
then $I (|f|) \leqslant r$, by positivity. By density of the Stone-Yosida
embedding, an integral extends uniquely to a positive linear functional on the
space of all continuous real-valued functions on the spectrum.

We present a geometric theory INT of integrals on $R$, much like the
description of Stone's maximal spectrum $\Max{R}$ in section~\ref{spectrum}.
In fact, the geometric theory Max will have one relation more than the theory
INT. This means that INT can be interpreted in Max, this interpretation
defines a frame map from INT to Max, and hence, a locale map from Max to INT.
The locale Max is a sublocale of INT. The inclusion is given by assigning to a
point its Dirac measure: $I_x (f) \assign f (x)$.

To wit, subbasic opens of INT, denoted by {\intl{p}{f}}, are indexed by
rational $p$ and $f$ in $R$. The set of its points will be $\compr{I}{p < I
(f)}$. Since $p < I (f)$ iff $0 < I (f - p)$, it is sufficient to treat basic
opens of the form $0 < I (f)$, written {\integral{f}}, where $P$ is a dummy
symbol. The points in this open are integrals $I$ such that $0 < I (f)$.

\begin{definition}
  The geometric theory INT is freely generated by symbols $\integral{f}$, $f$
  in $R$, and relations:

  \begin{descriptioncompact}
    \item[I1] $\integral{1} = \Top$;

    \item[I2] $\integral{f} \wedge \integral{\um f} = \Bot$;

    \item[I3] $\integral{f + g} \leqslant \integral{f} \vee \integral{g}$;

    \item[I4] $\integral{f} = \Bot$, if $f \leqslant 0$;

    \item[Cont] $P (f) = \bigvee_{s > 0} P (f - s)$.
  \end{descriptioncompact}
\end{definition}

\begin{lemma}
  \label{INT-rel}The relation $\integral{f} \leqslant \integral{g}$ if $f
  \leqslant g$ holds in INT.
\end{lemma}

\begin{proof}
  $P (f - g + g) \leqslant P (f - g) \vee P (g) = 0 \vee P (g)$.
\end{proof}

As before one proves that INT is compact completely regular by reducing
{\tmstrong{I1-4}}+{\tmstrong{Cont}} to {\tmstrong{I1-4}}. This result was
proved by Coquand~{\cite{Coquand:Stone}} who referred to the theory
{\tmstrong{I1-4}} as TOT, the theory of total orderings on an ordered vector
space. We have chosen the present presentation of INT since it makes compact
complete regularity easy to prove.

\begin{lemma}
  \label{model-integral}The theory INT is equivalent to the theory of
  normalized positive additive functionals:
  \begin{itemize}
    \item $I (f) \in \R $;

    \item $I (0) = 0$;

    \item $I (f + g) = I (f) + I (g)$;

    \item $I (f) \geqslant 0$ if $f \geqslant 0$;

    \item $I (1) = 1$.
  \end{itemize}
\end{lemma}

The notation above describes the locale with generators, $p < I (f)$ and $I
(f) < q$, for $f$ in $R$ and $p, q$ rational and $I$ is a dummy symbol, and
certain relations. For instance, the first axiom, $I (f) \in \R$, is a
shorthand for the relations:
\begin{itemize}
  \item $[p < I (f)] \leqslant [p' < I (f)]$ if $p' < p$;

  \item $[I (f) < q] \leqslant [I (f) < q']$ if $q < q'$;

  \item $p < I (f) = \bigvee_{p' > p} p' < I (f)$;

  \item $I (f) < q = \bigvee_{q' > q} q' < I (f)$;

  \item $1 = (p < I (f) \vee I (f) < q)$ if $p < q$;

  \item $0 = (q < I (f) \wedge I (f) < p)$ if $p < q$.
\end{itemize}
\begin{proof}
  We interpret $P (f)$ in INT as $I (f) > 0$ in the theory of positive
  additive functionals.

  For the converse, we define $p < I (f)$ as $\integral{f - p}$ and $I (f) <
  q$ as {\integral{q-f}}. Then
  \begin{enumerate}
    \item $- \varepsilon < I (f)$ if $f \geqslant 0$ and $\varepsilon > 0$.
    Proof: $\integral{f + \varepsilon - f} \leqslant \integral{f +
    \varepsilon} \vee \integral{- f}$ and $\integral{- f} = \Bot$.

    \item $\Top = s < I (f) \vee I (f) < t$, whenever $s < t$. Proof:
    $\integral{t - s} = \Top$.

    \item By Lemma~\ref{INT-rel}, if $s < I (f)$, then $t < I (f)$ for $t <
    s$. Similarly, if $I (f) < s$, then $I (f) < t$ for $s < t$.
  \end{enumerate}
  Combined with the continuity rule, this shows that $I (f)$ is a Dedekind
  cut.

  From {\tmstrong{I3}} we have $\integral{f} \leqslant \integral{\frac{1}{n}
  f}$. Hence, $1 - \frac{1}{n} < I (1) < 1 + \frac{1}{n}$, i.e. $I (1) = 1$.
  Similarly, $I (0) < \frac{1}{n}$.

  To prove additivity we combine Lemma~\ref{less-ineq} with {\tmstrong{I3}}
  and obtain $I (f + g) \leqslant I (f) + I (g)$. Conversely, the rule
  $\integral{f} \wedge \integral{g} \leqslant \integral{f + g}$ can be derived
  in INT: $f = f + g - g$, so $\integral{f} \leqslant \integral{f + g} \vee
  \integral{- g}$ and the result follows from $\integral{g} \wedge
  \integral{\um g} = \Bot$.
\end{proof}

Linearity readily follows from additivity, so the points of INT are integrals
and, conversely, every integral defines a point.

Usually, one proves that the space of integrals is compact by an appeal to the
Alaoglu theorem. Here we have shown that it is compact by construction. A
similar construction can be carried out for Alaoglu's theorem for compact
locales~{\cite{MulveyPelletier}}.

\subsection{Integrals on positive elements\label{subsec:intpos}}

Instead of starting with a positive linear functional, it will later be
convenient to work with its restriction to the positive elements.

\begin{lemma}
  An integral is fixed by its behavior on the positive elements. As such it is
  a function $I : R^+ \rightarrow \R^+$ such that $I (0) = 0$ and $I (f + g) =
  I (f) + I (g)$ and $I (1) = 1$.
\end{lemma}

The theory of these functionals is geometric, we call this theory INTPOS.

\begin{proposition}
  The geometric theories INT and INTPOS are biinterpretable.
\end{proposition}

\begin{proof}
  To obtain the integral from its positive part we define $I (f) \assign I
  (f^+) - I (f^-) .$
\end{proof}{\hspace*{\fill}}

\subsection{The space of valuations\label{subsec:val}}

\begin{definition}
  A {\tmem{valuation}} is a map $\mu : \tmop{Spec} (R) \rightarrow
  \R_{\tmop{low}}^+$ such that
  \begin{itemize}
    \item $\mu ( \Bot) = 0$, $\mu ( \Top) = 1$;

    \item $\mu (x) + \mu (y) = \mu (x \vee y) + \mu (x \wedge y)$ (the modular
    law);

    \item If $x \leqslant y$ in $\tmop{Spec} (R)$, then $\mu (x) \leqslant \mu
    (y)$ ($\mu$ is monotone);

    \item $\mu (D (a)) \leqslant \bigvee_{\varepsilon > 0} \mu (D (a -
    \varepsilon))$ {\tmem{Scott-continuous}}.
  \end{itemize}
\end{definition}

The theory of valuations is geometric, we call this theory VAL. Using
Lemma~\ref{less-ineq} we can formulate modularity in a way similar
to~{\cite{MoshierJung}}. We have defined the valuation only on the coherent
basis Spec(R) of $\Max{R}$, but it extends to the locale itself.
Alternatively, we could have used the same definition but with monotonicity
for the order on $\Max{R}$. This gives rise to the same locale of valuations:
If $D (a) \leqslant D (b)$ in $\Max{R}$, then $D (a - r) \leqslant D (b)$ in
$\tmop{Spec} (R)$ for all $r > 0$ and so $\mu (D (a - r)) \leqslant \mu (D
(b))$. By Scott-continuity we get $\mu (D (a)) \leqslant \mu (D (b))$.

This locale coincides with the locale of valuations on the locale $\Max{R}$ as
defined by Vickers~{\cite{Vickers:Integration}} for an arbitrary locale with
the difference that we require $\mu ( \Top) = 1$.
Vickers~{\cite[Prop.4.1]{Vickers:Integration}} already pointed out that we can
restrict to a base of the locale in order to obtain the locale of valuations
geometrically from (a presentation) of the locale.

Classically, the regular measures form a compact Hausdorff space. Hence,
classically, the locale of valuations on a compact completely regular locale
is again compact completely regular. The homeomorphism in the Riesz
representation theorem gives a constructive proof of this
fact.{\hspace*{\fill}}

\subsection{Statement of the theorem\label{statement}}

We are now ready to define the promised maps between integrals and valuations.
We give a syntactic bi-interpretation between two theories: the definition of
the maps will be geometric, but the reasoning that these maps actually satisfy
the required properties will be intuitionistic. For a general discussion of
such techniques see e.g.~{\cite[sec.4.5]{Vic:LocTopSp}}.

\subsubsection*{From integrals to valuations}

Given an integral on a Riesz space, we construct a valuation on the opens in
its spectrum:
\[ \mu_I (D (a)) \assign \sup \compr{I (n a^+ \wedge 1)}{n \in \N} \]
In section~\ref{sec:contmaps} we prove that this is well-defined, i.e.~that it
gives the same answer when $D (a) = D (b)$.

\subsubsection*{From valuations to integrals}

In order to define the converse interpretation we introduce some notations.
For $f$ in $R^+$ define the lower real $\Delta_f (r, s) \assign \mu (r < f <
s)$. Let $I = (r, s)$. Write $\Delta_f (I^c)$ for the lower real $\Delta_f (-
\infty, r) + \Delta_f (s, \infty)$ and $\Delta_f [I]$ for the {\tmem{upper}}
real $1 - \Delta_f (I^c)$.

The interpretation of INTPOS in VAL
\begin{eqnarray}
  I_{\mu} f & \assign & (\sup_{(s_i)} \sum s_i \Delta_f (s_i, s_{i + 1}),
  \inf_{(s_i)} \sum s_{i + 1} \Delta_f [s_i, s_{i + 1}]) \nonumber
\end{eqnarray}
the $(s_i)$ range over partitions over a fixed interval $[a, b]$ where $a < f
< b$. As is the case for $\mu_{\um}$ this is a disjunction over a concrete
countable set: a finite list of strictly increasing rationals.

Assuming the classical Riesz representation theorem it is easy to show that
these are indeed interpretations and that these maps are each other's inverses
as follows: For any $r > 0$ there is an $r$-approximation by sums $\sum s_i
\Delta_f (s_i, s_{i + 1})$ and $\sum s_{i + 1} \Delta_f [s_i, s_{i + 1}]$.
This follows from the usual classical proof of Riesz Theorem and the
possibility to choose $s_i$ as continuity points for the function
\[ s \mapsto \Delta_f (\um \infty, s) \]
By completeness of propositional
$\omega$-logic~{\cite{MakkaiReyes,ScottTarski}} and the validity of the
propositions in all models, i.e. measures or integrals, of the theory we see
that, classically, there should be a proof in the theory that these are indeed
interpretations. We will provide such a constructive proof in
Theorem~\ref{thm:Riesz}. This treatment is different from the classical one;
see e.g.~{\cite{Rudin}}. We take the topological/computational aspects into
account by distinguishing between lower reals and Dedekind reals, moreover we
do not use the extension of a valuation to a measure on the Borel sets. Our
result is more general: not only is it constructive, and hence valid in any
topos, but it also abstracts from a lattice of sets to a general
lattice.

\section{\label{Sec:Proof}Proof of the Riesz representation theorem}

\subsection{Formal simple functions\label{sub:simple1}}

We define formal simple functions on a distributive lattice $L$. All index
sets in this section are finite, i.e.~have a cardinality. We will use the
convention that a capital letter, say $I$, is a subset of the variables
indexed by the lower case letters, say $(x_i)$. For $(x_i)$ in $L$ we define
$x_I \assign \wedge \compr{x_i}{i \in I}$. Following Tarski~{\cite{Tarski}}
and Horn and Tarski~{\cite[Def 1.4]{HornTarski}} we define the free monoid $M
(L)$ such that the relation $x + y = x \vee y + x \wedge y$ holds. As Horn and
Tarski prove this is the monoid of formal sums $\sum x_i$, where $x_i$ in $L$,
with the following equality:

\begin{lemma}
  \label{simpl-eq}{\tmem{{\cite{HornTarski}}}} We have $\sum_{i \in I} a_i =
  \sum_{k \geqslant 1} \bigvee_{K \subset I, |K| = k} a_K .$
  Furthermore,$\sum_{i \in I} a_i = \sum_{j \in J} b_j$ iff $\bigvee_{K
  \subset I, |K| = k} a_K = \bigvee_{K \subset J, |K| = k} b_K$ for all $k
  \geqslant 1$.
\end{lemma}

\begin{definition}
  \label{def:monoid}Let $M (L)$ be the monoid of formal sums in $L$ modulo the
  relation $x + y = x \vee y + x \wedge y$. We define the pre-order
  \[ \sum x_i \leqslant \sum y_j \tmop{iff} \tmop{for} \tmop{all} I, x_I
     \leqslant \bigvee \compr{y_J}{|J| = |I|} . \]
  By Lemma~\ref{simpl-eq} $\leqslant$ is an order.
\end{definition}

The monoid $M (L)$ satisfies the cancellation property; see~{\cite{Tarski}}. \
For $k > 0$, $k x \leqslant 0$ iff $x = 0$. We add positive rational
coefficients {\emdash} that is, define a relation $\sum r_i x_i \leqslant \sum
s_j y_j$ {\emdash} by putting all the terms on one denominator. If $r$ in
$\Q^+$ and $x \leqslant y$, then $rx \leqslant ry$ and $x + z \leqslant y +
z$. When $L$ is a lattice of sets, this coincides with the usual ordering of
simple functions. We write $S^+ (L)$ for the {\tmem{positive simple
functions}} on $L$.

We write $r_I \assign \sum_{i \in I} r_i$. The following is direct.

\begin{lemma}
  $\sum r_i x_i \leqslant \sum s_j y_j$ iff for all $I$, $x_I \leqslant
  \bigvee_{J, r_I \leqslant s_J} y_J$.
\end{lemma}

\begin{lemma}
  The relation $\leqslant$ is transitive on $S^+ (L)$.
\end{lemma}

\begin{proof}
  Suppose that $\sum r_i a_i \leqslant \sum s_j b_j \leqslant \sum t_k c_k$.
  By Lemma~\ref{simpl-eq}, for all $I$, $a_I \leqslant \bigvee_{J, r_I
  \leqslant s_J} b_J$ and for all $J$, $b_J \leqslant \bigvee_{K, s_J
  \leqslant t_K} c_K$. So, $a_I \leqslant \bigvee_{J, K, r_I \leqslant s_J,
  s_J \leqslant t_K} c_K$.
\end{proof}

\subsection{Extending valuations to simple functions\label{Val-simple}}

We now consider the case where $L$ is $\tmop{Spec} (R)$. We extend $\mu$ to an
additive functional from the formal sums to the lower reals. This extension
satisfies the modular law and hence so does the extension to the simple
functions:

\begin{lemma}
  If $\sum r_i x_i \leqslant \sum s_j y_j$, then $\mu ( \sum r_i x_i)
  \leqslant \mu ( \sum s_j y_j)$. So, $\mu$ is well-defined on $S^+ (L)$: if
  $k = l$, then $\mu (k) = \mu (l)$.
\end{lemma}

\begin{proof}
  By bringing all the terms on one denominator we can dispose of all the
  scalars. Hence our goal will be to prove: If $\sum x_i \leqslant \sum y_j$,
  then $\mu ( \sum x_i) \leqslant \mu ( \sum y_j)$. To see this we have
  \[\mu ( \sum x_i)\,\, =\,\,  \mu ( \sum_{k \geqslant 1} \bigvee_{|K| = k} x_K)
  \,\,\leqslant\,\,  \mu ( \sum_{k \geqslant 1} \bigvee_{|K| = k} y_K)
  \,\,=\,\,\mu ( \sum y_j).\proved\]
   %
\end{proof}

Consider the dual lattice $L'$ of $\tmop{Spec} (R)$. We define $\mu (\neg x)$
as the {\tmem{upper}} real $1 - \mu (x)$. This definition is naturally
extended to the formal simple functions $S^+ (L')$: $\mu ( \sum s_i (\neg
x_i)) = ( \sum s_i) - \mu ( \sum s_i x_i)$. However, we will not be able to
define the valuation of a sum of mixed open and closed elements.

\subsection{Simple functions on the spectrum of a Riesz
space\label{subsec:simple}}

We now consider the case where $L$ is the Boolean algebra freely generated by
$\tmop{Spec} (R)$. Let $f$ be in $R$. We denote the open $D (f - r)$ by $(f >
r)$ and $D (r - f)$ by $(f < r)$ and the complement of $(f > r)$ by $(f
\leqslant r$) and the complement of $(f < r)$ by $(f \geqslant r$).

We want to express the pointwise order relation between a positive simple
function and a positive element of the Riesz space considered as continuous
functions on the spectrum $\Max{R}$. However, for the sake of geometricity, we
use the order of $\tmop{Spec} (R)$ instead. Hence we are working with a
coherent approximation to the pointwise order.

We define the relation $\sum r_i x_i \leqslant f$ as: for all $I$, $x_I
\leqslant (r_I \leqslant f)$ and the relation $f \leqslant \sum s_j y_j$ as:
$\Top = \bigvee_J ((f \leqslant s_J) \wedge y_J)$.

\begin{lemma}
  \label{D(a)}If $a \leqslant 1$, then $a \leqslant D (a)$.
\end{lemma}

\begin{proof}
  We need to prove that $(a \leqslant 0) \vee ((a \leqslant 1) \wedge D (a)) =
  1$. We simplify this statement:

  $(a > 0) \leqslant (a \leqslant 1) \wedge D (a)$

  $(a > 0) \leqslant (a \leqslant 1)$

  $(a > 0) \wedge (a > 1) = 0$

  The last statement follows from the hypothesis $a \leqslant 1$.
\end{proof}

When $\Max{R}$ is spatial, as is the case in the presence of the axiom of
choice, by Stone-Yosida, $f$ may be interpreted as a continuous function on
$\Max{R}$ and the order above corresponds to a {\tmem{coherent approximation}}
of the pointwise ordering of functions when the simple function $\sum r_i x_i$
is interpreted as the linear combination of the characteristic functions
associated to the sets $x_i$.

\begin{lemma}
  Suppose that $\sum r_i x_i \leqslant \sum s_j y_j$ and $\sum s_j y_j
  \leqslant f$. Then $\sum r_i x_i \leqslant f$.
\end{lemma}

\begin{proof}
  We have $x_I \leqslant \bigvee_{J, r_I \leqslant s_J} y_J$ and $y_J
  \leqslant (s_J \leqslant f)$. So
  \[ x_I \leqslant \bigvee_{J, r_I \leqslant s_J} (s_J \leqslant f) \leqslant
     \bigvee_{J, r_I \leqslant s_J} (r_I \leqslant f) \leqslant (r_I \leqslant
     f) . \proved\]
\end{proof}

\begin{lemma}
  Suppose that $f \leqslant \sum r_i x_i$ and $\sum r_i x_i \leqslant \sum s_j
  y_j$. Then $f \leqslant \sum s_j y_j$.
\end{lemma}

\begin{proof}
  We have $\Top = \bigvee_I ((f \leqslant r_I) \wedge x_I)$ and $x_I \leqslant
  \bigvee_{J, r_I \leqslant s_J} y_J$. So
  \begin{eqnarray}
   \Top = \bigvee_I ((f \leqslant r_I) \wedge x_I)&=& \bigvee_I ((f \leqslant
    r_I) \wedge \bigvee_{J (I), r_I \leqslant s_{J (I)}} y_{J (I)}) \nonumber \\
    &    \leqslant & \bigvee_I \bigvee_{J (I), r_I \leqslant s_{J (I)}} (f
    \leqslant s_{J (I)}) \wedge y_{J (I)}
     \leqslant  \bigvee_J ((f \leqslant s_J) \wedge y_J) . \nonumber
  \end{eqnarray}
\end{proof}

It is clear that if $\sum r_i x_i \leqslant f \leqslant g$, then $\sum r_i x_i
\leqslant g$, and if $f \leqslant g \leqslant \sum r_i x_i$, then $f \leqslant
\sum r_i x_i$.

\begin{lemma}
  \label{lem:ineq-simpl}If $\sum r_i x_i \leqslant f \leqslant \sum s_j y_j$,
  then $\sum r_i x_i \leqslant \sum s_j y_j$.
\end{lemma}

\begin{proof}
  We have for all $I$, $x_I \leqslant (r_I \leqslant f)$ and $\Top = \bigvee_J
  ((f \leqslant s_J) \wedge y_J)$. Then
  \[ x_I \leqslant (r_I \leqslant f) = (r_I \leqslant f) \wedge \bigvee_J ((f
     \leqslant s_J) \wedge y_J) \leqslant \bigvee_{J, r_I \leqslant s_J} y_J .\proved
  \]

\end{proof}

\begin{lemma}
  \label{partition}Let $0 \leqslant f \leqslant b$ and let $s_i$ be a
  partition of $[0, b]$. Then
  \[ \sum s_i (s_i < f < s_{i + 1}) \leqslant f \leqslant \sum s_{i + 1} (s_i
     \leqslant f \leqslant s_{i + 1}) . \]
\end{lemma}

\begin{proof}
  To prove the first inequality, we write $x_i \assign_{} (s_i < f < s_{i +
  1})$. The $x_i$ are disjoint and
  \[ x_i \leqslant (s_i < f) \leqslant (s_i \leqslant f) . \]
  To prove the second inequality we write $t_i \assign (s_i \leqslant f
  \leqslant s_{i + 1})$. Then $(f \leqslant s_{i + 1}) \wedge y_i = y_i$ and
  $1 = \bigvee y_i$, since $s_0, ..., s_n$ is a partition of $[0, b]$.
\end{proof}

The following results have a direct proof.

\begin{lemma}
  \label{sum}If $l_1 \leqslant f_1$ and $l_2 \leqslant f_2$, then $l_1 + l_2
\leqslant f_1
  + f_2$. Similarly, if $f_1 \leqslant k_1$ and $f_2 \leqslant k_2$, then $f_1 +
f_2 \leqslant k_1 + k_2$.
\end{lemma}

%

The spectrum of a Riesz space is completely regular as the following simple
formulation of the Urysohn's Lemma shows.

\begin{lemma}
  \label{Ury}Let $a$ in $R^+$ and $\varepsilon > 0$. Then $D (a - \varepsilon)
  \leqslant \frac{1}{\varepsilon} (a \wedge \varepsilon) \leqslant D (a)$.
\end{lemma}

\begin{proof}
  For the first inequality we need to prove that $D (a - \varepsilon)
  \leqslant (1 \leqslant \frac{1}{\varepsilon} (a \wedge \varepsilon))$. Since
  the right hand side is a formal complement this means, $D (a - \varepsilon)
  \wedge (1 > \frac{1}{\varepsilon} (a \wedge \varepsilon)) = 0.$ Now, $(1 >
  \frac{1}{\varepsilon} (a \wedge \varepsilon)) = (\varepsilon > a) = D
  (\varepsilon - a)$.

  The second inequality follows from Lemma~\ref{D(a)}: $\frac{1}{\varepsilon}
  (a \wedge \varepsilon) \leqslant D ( \frac{1}{\varepsilon} (a \wedge
  \varepsilon)) = D (a)$.
\end{proof}

\subsection{\label{sec:Delta}$\Delta$-functions}

In this subsection we fix $f \geqslant 0$ in $R$ and a valuation $\mu$. We
define the lower real $\Delta (r, s) = \mu (r < f < s)$ and the upper real
$\Delta [r, s] = 1 - \Delta (- \infty, r) - \Delta (s, \infty)$ as in
Section~\ref{statement}. Intuitively, the function $\Delta$ represents the
function $\alpha (s) = \mu (f < s)$ which is used in the definition of the
integral as a Stieltjes integral $\int f \mathd \mu = \int s \mathd \alpha
(s)$. The functions $\Delta$ satisfies:
\begin{enumerate}
  \item $\Delta (0, b) = 1$ for some $b$;

  \item $\Delta (r, s) \leqslant 1$;

  \item $\Delta (r, s) \geqslant 0$;

  \item $\Delta (r, s) + \Delta (s, t) = \Delta (r, t) - \Delta [s] ;$

  \item $\Delta (r', s') \leqslant \Delta (r, s)$ whenever $r \leqslant r' <
  s' \leqslant s$;

  \item $\Delta (r, s') + \Delta (r', s) = \Delta (r, s) + \Delta (r', s')$
  whenever $r < r' < s' < s$;

  \item $\Delta (r, s) = \bigvee \compr{\Delta (r', s')}{r < r' < s' < s} .$
\end{enumerate}
In 4, $\Delta [s] \assign \Delta [s, s]$.

We write $(r', s') \ll (r, s)$ for $\text{$r < r' < s' < s$}$. As before, we
write $\Delta (I)$ for $\Delta (r, s)$, if $I = (r, s)$.

\begin{lemma}
  If $I \ll J$ and $p < q$, then $\Delta (J) > p$ or $\Delta [I] < q$.
\end{lemma}

\begin{proof}
  Since $\Delta (I^c) + \Delta (J) \geqslant 1 > p + (1 - q)$.
\end{proof}

We now prove `a non-increasing function is continuous in a dense set of
points' in a pointfree way.

\begin{theorem}
  \label{thm:Delta}Let $N \in \N$ and $I = (r, s)$ be an open interval. Then
  there exists an interval $J \ll I$ such that $\Delta [J] < \frac{1}{N}$.
\end{theorem}

\begin{proof}
  Choose $2 N$ disjoint intervals $I_i$ in $I$ and choose $2 N$ intervals $J_i
  \ll I_i$. For each $i$, $\Delta (I_i) > \frac{1}{2 N}$ or $\Delta [J_i] <
  \frac{1}{N}$. It is impossible that the former case occurs all the time,
  therefore the latter case occurs at least once.
\end{proof}

It follows classically that $\mu (r < f \leqslant s) = \inf_{s' > s} \mu (r <
f < s')$. The approximations to this infimum are explicit in the following
proposition which assigns a Dedekind real to $\Delta$. The interpretation of
this real is the Stieltjes integral$\int s \mathd \alpha (s)$, where $\alpha$
is a non-decreasing function connected to $\Delta$.

\begin{proposition}
  \label{prop:sup-exists}The pair
  \[ ( \compr{p}{p < \sum s_i \Delta (s_i, s_{i + 1})}, \compr{q}{\sum s_{i +
     1} \Delta [s_i, s_{i + 1}] < q}), \]
  where $s_i$ ranges over finite partitions of $[0, b]$, defines a Dedekind
  real.
\end{proposition}

\begin{proof}
  We first prove that the upper and lower cut come arbitrary close: There
  exists $(s_i)$ such that $\sum s_{i + 1} \Delta [s_i, s_{i + 1}] - \sum s_i
  \Delta (s_i, s_{i + 1})$ is small. To wit, given $\varepsilon > 0$, use
  Theorem~\ref{thm:Delta} to choose a partition $s_i$ of $[a, b]$ such that
  $|s_{i + 1} - s_i | < \varepsilon$ and $\sum \Delta [s_i] < \varepsilon$.
  Then
  \begin{eqnarray}
    \sum s_{i + 1} \Delta [s_i, s_{i + 1}]\!\!\! &-&\!\!\! \sum s_i \Delta (s_i, s_{i + 1})\nonumber \\
    & \leqslant & \sum (s_{i + 1} - s_i) \Delta [s_i, s_{i + 1}] + \sum s_i
    (\Delta [s_i] + \Delta [s_{i + 1}]) \nonumber\\
    & \leqslant & \varepsilon \sum \Delta [s_i, s_{i + 1}] + 2 b \varepsilon
    \nonumber\\
    & \leqslant & \varepsilon (1 + \sum \Delta [s_i]) + 2 b \varepsilon
    \leqslant \varepsilon (1 + \varepsilon) + 2 b \varepsilon . \nonumber
  \end{eqnarray}
  We now prove that the lower cut is below the upper cut. By
  Lemma~\ref{partition},
  \[ l \assign \sum s_i (s_i < f < s_{i + 1}) \leqslant f \leqslant \sum s_{i
     + 1} (s_i \leqslant f \leqslant s_{i + 1}) = : k. \]
  Write $y_i \assign (f < s_i) \vee (s_{i + 1} < f)$. By
  Lemma~\ref{lem:ineq-simpl}, $l + \sum s_{i + 1} y_j \leqslant \sum s_{i +
  1}$, so $\mu (l) + \mu ( \sum s_{i + 1} y_j) \leqslant \mu ( \sum s_{i +
  1})$. The conclusion, $\mu (l) \leqslant \mu (k)$, follows from
  Lemma~\ref{low-up}.
\end{proof}

The previous proposition contains the essence of Bishop's profile theorem;
see~{\cite{Bishop/Bridges:1985}}. It is the crucial step in the proof that
$I_{\mu}$ is a function; see Lemma~\ref{lem:VI}.

\subsection{\label{sec:contmaps}Continuous maps}

We are now ready to show that the maps $\mu_I$ and $I_{\mu}$ defined above
indeed map integrals to valuations, and vice versa. We need to check that the
interpretations of all the axioms hold.

We first repeat the definition from section~\ref{statement}:
\[ \mu_I (D (a)) \assign \sup \compr{I (n a^+ \wedge 1)}{n \in \N} . \]
We observe that $f \leqslant D (a)$ iff there exists $n$ such that $f
\leqslant n a^+ \wedge 1$. Hence, $\sup \compr{I (n a^+ \wedge 1)}{n \in \N} =
\sup \compr{I (f)}{f \leqslant D (a)} .$

The map $\mu_I$ extends to the positive simple functions:
\[ \mu_I ( \sum r_i D (a_i)) = \sup \compr{I ( \sum r_i (n a^+_i \wedge 1))}{n
   \in \N} . \]
\begin{lemma}
  $\mu_I$ is a valuation.
\end{lemma}

\begin{proof}
  To prove modularity we observe that
  \begin{gather*}
  (n a \wedge 1) + (n b \wedge 1) = (n (a \wedge b) \wedge 1) + (n (a \vee b) \wedge 1)\\
  \tag*{\text{and hence}}\\
  I (n a \wedge 1) + I (n b \wedge 1) = I (n (a \wedge b) \wedge 1) + I (n
     (a \vee b) \wedge 1) .
  \end{gather*}
  For monotonicity: If $f \leqslant x$ and $x \leqslant y$, then $f \leqslant
  y$. Finally, regularity, $\mu (D (a)) = \sup_{r > 0} \mu (D (a - r))$ is
  direct.
\end{proof}

We generalize the definition of $I_{\mu}$ in section~\ref{statement} to
arbitrary simple positive functions:
\begin{eqnarray}
  I_{\mu} (f) & = & (\sup \compr{\mu (l)}{l \leqslant f, l \in S^+ (L)}, \inf
  \compr{\mu (k)}{f \leqslant k, k \in S^+ (L')}) . \nonumber
\end{eqnarray}
We will prove that the supremum and the infimum over the restricted sets of
simple functions used in section~\ref{statement} already form a Dedekind real
and hence the two definitions coincide.

\begin{lemma}
  \label{lem:VI}$I_{\mu}$ is an integral.
\end{lemma}

\begin{proof}
  To prove that $I$ maps to the Dedekind reals: Let $f \in R^+$ and choose $b
  \geqslant f$. By Proposition~\ref{prop:sup-exists} $(\sup \sum s_i \Delta
  (s_i, s_{i + 1}), \inf \sum s_{i + 1} \Delta [s_i, s_{i + 1}])$ is a
  Dedekind real: the lower cut is below the upper cut and both cuts `kiss'.

  To prove additivity, by Lemma~\ref{sum}, if $l_1 \leqslant f$ and $l_2
  \leqslant g$, then $l_1 + l_2 \leqslant f + g$. Hence $I (f) + I (g)
  \leqslant I (f + g)$. Conversely, if $f \leqslant k$ and $g \leqslant l$,
  then $f + g \leqslant k + l$ \ and hence $I (f + g) \leqslant I (f) + I
  (g)$.
\end{proof}

\subsection{Homeomorphism}

We prove that there is a homeomorphism between the integrals on a Riesz space
and the valuations on the opens of the spectrum.

\begin{theorem}
  \label{thm:Riesz}{\tmem{[Riesz representation theorem]}} Let $R$ be a Riesz
  space with a strong unit. The theory of valuations on its spectrum is
  equivalent to the theory of integrals on $R$. It follows that the
  corresponding compact completely regular locales are homeomorphic.
\end{theorem}

\begin{proof}
  That is, we claim that $I_{\mu_J} = J$ and $\mu_{I_{\nu}} = \nu$.
  \begin{eqnarray}
    I_{\mu_J} (f) & = & \sup ( \compr{\mu_J (l)}{l \leqslant f})
    \nonumber\\ & = &
     \sup \compr{J (g)}{g \leqslant l \leqslant f} \nonumber\\
    & \geqslant & \sup \compr{J (f - \varepsilon)}{\varepsilon > 0} = J (f)
    \nonumber
  \end{eqnarray}
  For the inequality we observe that for each $\varepsilon > 0$, $f -
  \varepsilon \leqslant \sum_{n \geqslant 1} \frac{\varepsilon}{2} ((n +
  \frac{1}{2}) \frac{\varepsilon}{2} < f) \leqslant f$.

  The other inequality is trivial.

  Conversely,
  \begin{eqnarray}
    \mu_{I_{\nu}} (k) & = & \sup \compr{I_{\nu} (f)}{f \leqslant k}
    \nonumber\\& = &
    \sup \compr{\nu (l)}{l \leqslant f \leqslant k} \nonumber\\
    & \geqslant & \sup \compr{\nu (l)}{l \ll k} = \nu (k) \nonumber
  \end{eqnarray}
  Where $l \ll k$ means $k = \sum s_i D (a_i)$ and $l = \sum s_i D (a_i -
  \varepsilon)$. By the Urysohn Lemma~\ref{Ury} there exists $f$ in $R$ such
  that $l \leqslant f \leqslant k$. The other inequality is trivial.
\end{proof}

\section{Related work}

Vickers~{\cite{Vickers:Integration}} presents another variant of the Riesz
representation theorem. His construction works for locales which are not
necessarily compact completely regular. However, his integrals have their
values in the lower (or upper) reals, as opposed to the Dedekind reals. A
locale of valuations was first presented by Heckman~{\cite{Heckmann}}.

The present homeomorphism has already been applied in a
{\tmem{non-commutative}} context of quantum theory~{\cite{HLS}} where it
provides an isomorphism between quasi-states and certain valuations.
Quasi-states are used in the algebraic foundations of quantum mechanics.

\section{Conclusions}

The present construction was motivated by Bishop's bijection between measures
and integrals~{\cite{Bishop/Bridges:1985}}. Bishop's forces the measure of a
measurable set to be a Dedekind real. This is somewhat inconvenient in
practice since for a measurable function $f$ the measure of $[f \geqslant s]$
need not be Dedekind in general. We believe that the present theory allows for
a smoother development of, at least, the abstract functional analytic aspects
of Bishop's measure theory.

\section{Acknowledgements}

We would like to thank Alex Simpson, Steve Vickers and the referees for
comments on the presentation of this paper.

\bibliographystyle{jloganal}

\bibliography{profile}

\end{document}